\documentclass[12pt]{amsart} 
\usepackage{amsfonts} 
\usepackage{graphicx}
\usepackage{color}
\usepackage{wrapfig}
\usepackage{hyperref}
\usepackage{todonotes}
\hypersetup{backref,colorlinks=true,citecolor=blue,linkcolor=blue}

\usepackage{mdframed}

\theoremstyle{definition}

\newcommand{\nc}{\newcommand}

\nc{\p}[1]{\medskip\noindent{\em #1.}}
\nc{\margin}[1]{\marginpar{\scriptsize #1}}

\begin{document}

\title{Alternating Knots}

\author{William W. Menasco}
\address{Department of Mathematics \\ University at Buffalo---SUNY}
\email{menasco@buffalo.edu}

\begin{abstract}
This is a  short expository article on alternating knots and is to appear in the Concise Encyclopedia of Knot Theory.
\end{abstract}

\maketitle

\setcounter{enumi}{0}
\section*{Introduction}
\label{intro}

\begin{figure}[h]
\begin{center}
{\includegraphics[width=8.8cm,height=5.6cm]{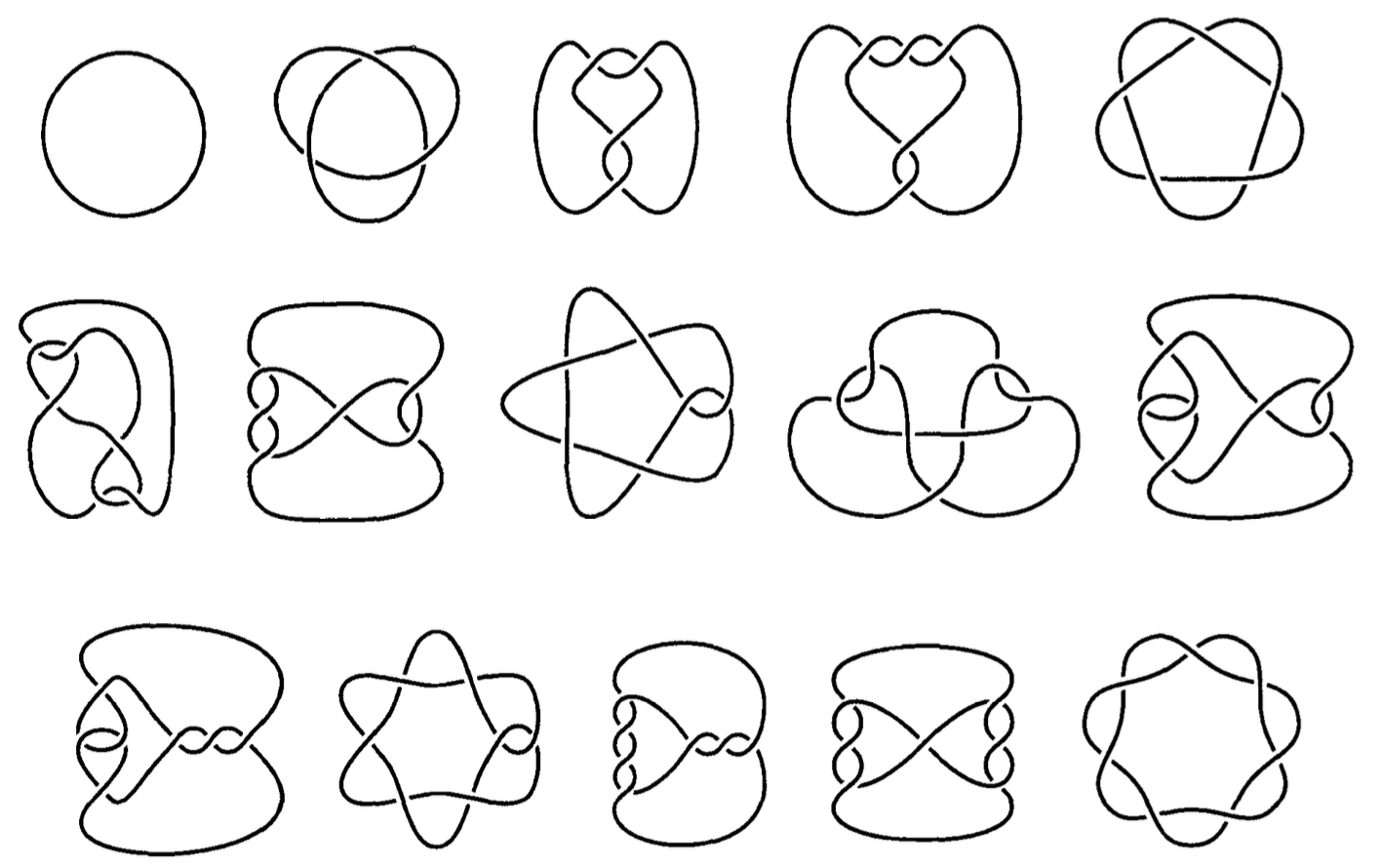}}
\end{center}
\caption{P.G. Tait's first knot table where he lists all knot types up to $7$ crossings.  (From reference [6], courtesy of J. Hoste, M. Thistlethwaite and J. Weeks.)}
\label{Table}
\end{figure}

A knot $K \subset S^3$ is \emph{alternating} if it has a regular planar diagram $D_K \subset \mathbb{P}(\cong S^2) \subset S^3$
such that, when traveling around $K$ , 
the crossings alternate, over-under-over-under, all the way along $K$ in $D_K$.  Figure \ref{Table} show the first $15$ knot types in P. G. Tait's earliest table
and each diagram exhibits this alternating pattern.
This simple definition is very unsatisfying.
\emph{A knot is alternating if we can draw it as an alternating diagram?}  There is no mention of any geometric structure. 
Dissatisfied with this characterization of an alternating knot, Ralph Fox (1913-1973) asked: "What is an alternating knot?"

\begin{figure}[h]
    \includegraphics[height=2cm]{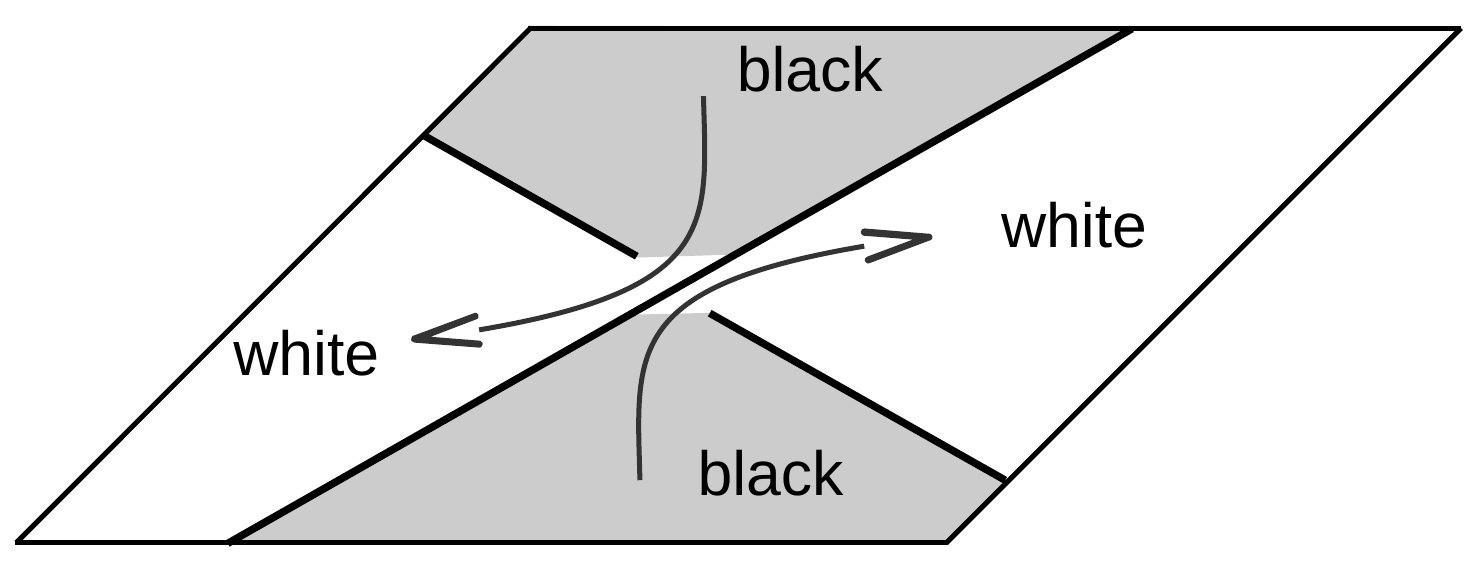}
  \caption{Going from a black to white region near a crossing.}
  \label{Crossing}
\end{figure}

Let's make an initial attempt to address this dissatisfaction by giving a different characterization of an alternating diagram that is
immediate from the over-under-over-under characterization.
As with all regular planar diagrams of knots in $S^3$, the regions of
an alternating diagram can be colored in a checkerboard fashion.  Thus, at each crossing (see figure \ref{Crossing}) we will have
``two'' white regions and ``two'' black regions coming together with similarly colored regions being kitty-corner to each other.
\begin{figure}[h]
    \includegraphics[height=2cm]{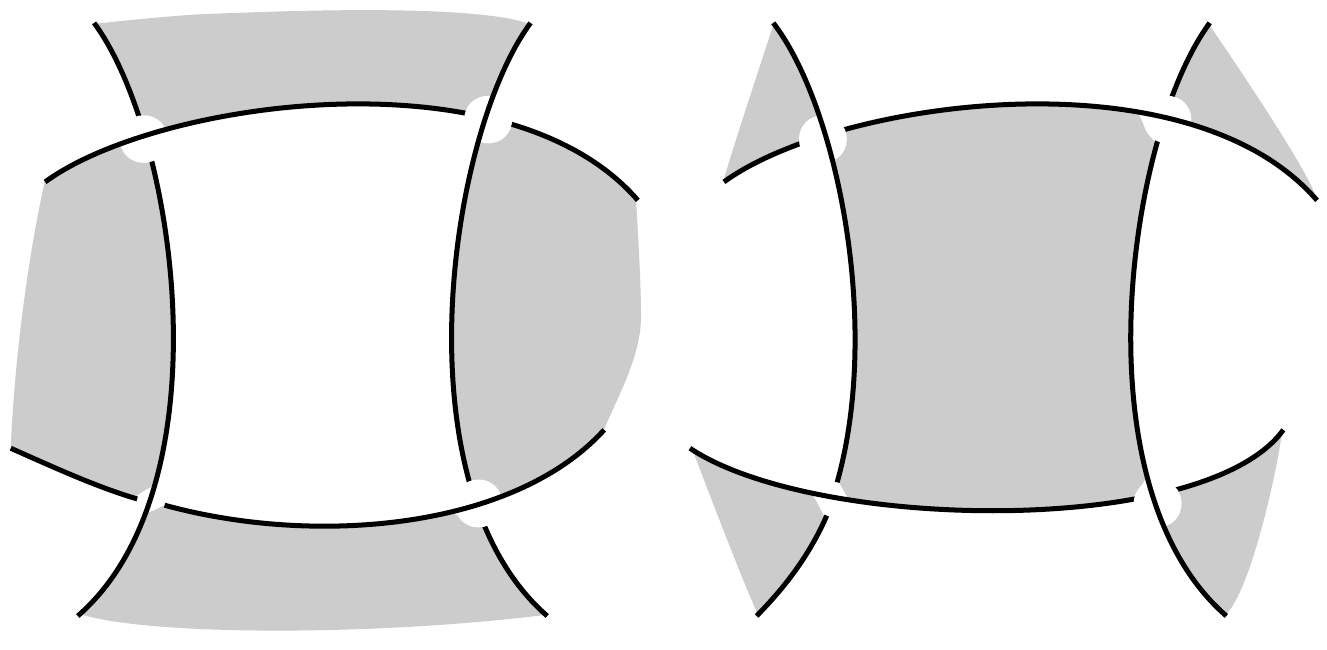}
 \caption{Black and white regions.}
   \label{regions}
\end{figure}
(I will explain my use of quotation marks momentarily.)  In the local picture of figure \ref{Crossing}, we will make the convention that
when standing in a black region and walking to a white region by stepping over the under-strand of the knot the over-strand
is to our left.  (Notice this is independent of which side of $\mathbb{P}$ we stand.)
Now the feature for an alternating diagram is, if one crossing of a black region satisfies this
local scheme then every crossing of that black region satisfies this scheme.  Since both black regions at a crossing satisfy our scheme
then this scheme is transmitted to every black region of an alternating diagram.  To summarize, a diagram of a knot is \emph{alternating}
if we can checkerboard color the regions such that when we stand in any black region, we can walk to an adjacent white
region by stepping over the under-strand of a crossing while having the over-strand of that crossing positioned to our left.
From the perspective of walking from a white regions to a black regions, the over-crossing is to our right.  Again, one
should observe that this coloring scheme is consistent independent of which side of $\mathbb{P}$ one stands.

\begin{figure}[h]
    \includegraphics[height=.9cm]{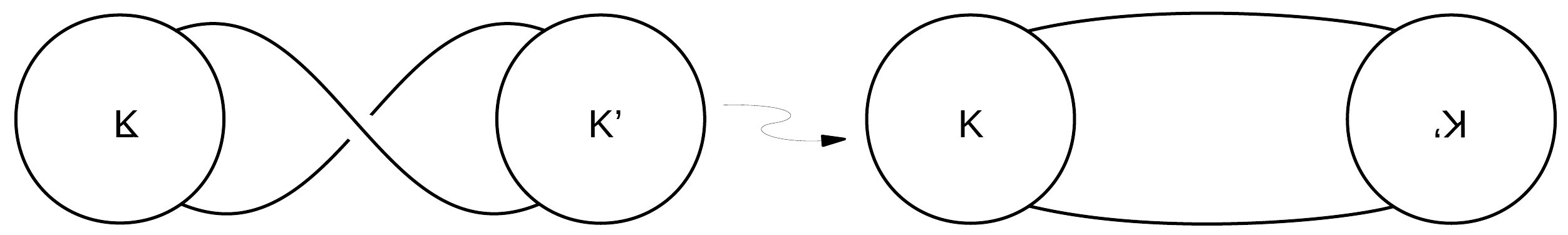}
  \caption{Eliminating a nugatory crossing.}
  \label{Nugatory}
\end{figure}

Coming back now to our previous use of quotation marks---``two'' white regions and ``two'' black regions---we must allow for the possibility
that ``two'' kitty-corner regions may be the same region as shown in figure \ref{Nugatory}.  This possibility forces the occurrence of a \emph{nugatory}
crossing in the diagram of $D_K$.  Such a crossing can be eliminated by a $\pi$-rotation of $K^\prime$ for half the diagram.
Observe that such
a $\pi$-rotation will take an alternating diagram to an alternating diagram with one less crossing.  So the elimination of nugatory crossings
must terminate and we call a diagram with no nugatory crossings \emph{reduced}.  Thus, our study of alternating diagrams is legitimately
restricted to ones that are reduced.  And, let's add one further restriction, that a diagram $D_K$ be connected, since one that is not connected is
obviously representative of a split link.
This requirement implies that all regions of a diagram $D_K$ are discs

With all of the above introductory material in place, and restricting to reduced connected alternating knot diagrams,
we now give in to the urge that anyone who has drawn a M\"obius band will have.  At each
crossing we place a half-twisted white band connecting the two white regions; similarly, we place a half-twisted black band
connecting the two black regions.  This construction gives us a white surface, $\omega_K$, and a black surface, $\beta_K$, each one having our knot $K$ as
its boundary---they are \emph{spanning surfaces} of $K$.
The white and black bands of each crossing intersect in a single arc that has one endpoint on the under-crossing
strand and its other endpoint on the over-crossing strand.  Thus, we can reinterpret a checkerboard $D_K$ diagram to be that of a black and white
surface with each crossing concealing an intersection arc since our eye's viewpoint would be one where we were looking straight down each arc.

Here is where I really wish to start our story of alternating knots.  (The story can be expanded to alternating links, but we will leave that for another time.)
At the center of this story lies three mysteries, the three conjectures
of the Scottish physicist Peter G. Tait (1831-1901).
As with all good conjectures, Tait's three point to deeper mathematics.  Alternating knots have proven particularly
well behaved with respect to tabulation, classification, and computing and topologically interpreting algebraic invariants.  
Moreover, for hyperbolic alternating knots their alternating diagrams are more
closely tied to their hyperbolic geometric structure than other knots.  As a collection, alternating knots have supplied researchers with a ready population for
experimentation and testing conjectures.  And, our two colored surfaces, $\beta_K$ and  $\omega_K$, are rightfully thought of as the heroes
the story for they are at the center of establishing the validity of the Tait conjectures and, finally, giving a satisfactory answer to Fox's question.

\section*{The Tait conjectures}\label{TC}

The early efforts of knot tabulation came initially from Tait.  He was motivated by Lord Kelvin's
program for understanding the different chemical elements as different knotted vortices in ether.
Without the aid of any theorems from topology, Tait published in 1876 his first knot table which containing the $15$ knot types through seven crossings in figure \ref{Table}.
Specifically, Tait enumerated all possible diagrams up to a seven crossings and then grouped together those diagrams that represented the same knot type.  For example,
the right-handed trefoil has two diagrams, one as a closed $2$-braid and one as a $2$-bridge knot.  Similarly, the figure-eight two diagrams, one coming from
a closed $3$-braid and one coming from a $2$-bridge presentation.  His grouping of these $15$ knot types is consistent with todays modern tables.
However,
his grouping did contained more than these $15$ since he did not possess the notion of prime/composite knots---he included the composite sums
of the trefoil with itself and the figure-eight.

\begin{figure}[h]
    \includegraphics[height=2cm]{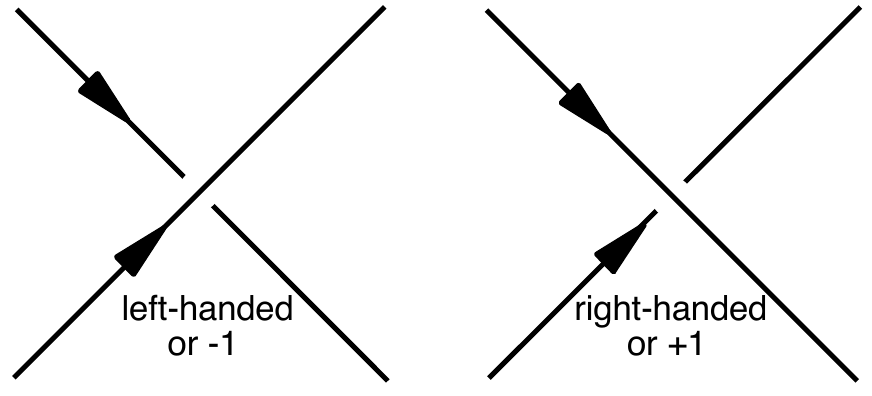}
  \caption{Left and right hand crossings.}
  \label{parity}
\end{figure}

Tait also experimented with the assigning of orientations to the knot diagrams---giving a defined direction to the knot along which to travel---which
induced a handedness to each crossing.  Thus, for right-handed crossings one
could assign a $+1$ parity and a $-1$ parity for left-handed crossings.  (See figure \ref{parity}.)  The sum of these parities was the \emph{writhe}  of a diagram.
For example, independent of orientation the trefoil in Tait's table will have three right-handed crossings yielding a writhe of $+3$,
whereas the figure-eight will always have two right-handed and two left-handed crossings for a writhe of $0$.

\begin{figure}[h]
    \includegraphics[height=.8cm]{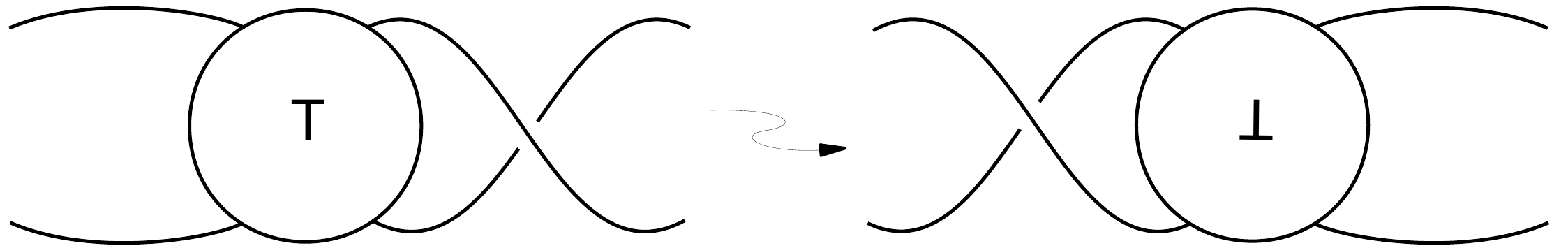}
  \caption{A flype moves the crossing from right to left.}
  \label{Flype}
\end{figure}

Tait also observed that one could reduce crossing number by the elimination of nugatory crossings.
Finally, Tait discovered an operation on an alternating diagram---a \emph{flype}---that preserved the number of crossings, the writhe and, more importantly, the knot type.
(See figure \ref{Flype}.)

From the diagrams and their groupings into just $15$ knot types, Tait had amassed a sizable about of data.
And, from this data Tait 
proposed a set of conjectures:
\vspace{10pt}

\begin{mdframed}[backgroundcolor=gray!20]
\centerline{\bf The Tait Conjectures}
\begin{enumerate} 
\item[\rm (T1)]  A reduced alternating diagram has minimal crossing number.
\item[\rm (T2)] Any two reduced alternating diagrams of the same knot type have the same writhe.
\item[\rm (T3)] Any two reduced alternating diagrams of a given knot type are related by a sequence of flypes.
\end{enumerate} 
\end{mdframed}
\vspace{10pt}

One should observe that T1 and T3 imply T2.  The solutions to these three conjectures would have to wait until after Vaughan Jones' work on his knot polynomial in 1984.

\section*{Surfaces and alternating knots}
\label{surfaces}
During the $20^{\rm th}$ Century, work on understanding the topology and geometry of alternating knots remained active.  But, for our
story line we start in the 1980's and describe how our colored surfaces, $ \beta_K$ and $\omega_K$
control the behavior of closed embedded surfaces in an alternating knot complement $S^3 - K$.  Any embedded closed surface, $\Sigma$,
in $S^3$ is necessarily orientable and
the ones that capture meaningful information about the topology of a knot complement are \emph{incompressible}.  That is,
every simple closed curve (s.c.c.) $c \subset \Sigma (\subset S^3 - K)$ which bounds an embedded disc in $S^3 - K$, also bounds a sub-disc of
$\Sigma$.  To this end, it is worth observing from the start that
$S^3 - (\beta_K \cup \omega_K)$ is the disjoint union of two open $3$-balls which we denote by $B^3_N$ ($N$ for north) and
$B^3_S$ ($S$ for south).

Given an incompressible surface $\Sigma \subset S^3 - K$ our first task is to put $\Sigma$ into normal position with respect to $\beta_K$
and $\omega_K$.  Using some basic general position arguments for incompressibility, this means that we can assume:
(1) $\beta_K \cap \Sigma$ and $\omega_K \cap \Sigma$ are collections of s.c.c.'s; (2) any curve from either collection cannot be totally contained
in a single black or white region used to construct $\beta_K$ or $\omega_K$---it must intersect some of the half-twisted bands of our colored surfaces---and;
(3) $B^3_N \cap \Sigma$ and $B^3_S \cap \Sigma$ are collections of open $2$-discs which we will call \emph{domes}.

Of particular interest is the the behavior of the boundary curves of the domes.  Let $\delta \subset B^3_{N} \cap \Sigma$ be a dome and consider the
curve $\partial (\bar{\delta}) \subset \beta_K \cup \omega_K$.  When viewed from inside $B^3_{N}$, the s.c.c. $\partial (\bar{\delta})$ will be seen as a union
of arcs---an arc in a white region adjoined to an arc in a black region adjoined to an arc in a white region, etc.---where the adjoining of two
consecutive arcs occurs within an intersection arc of $\beta_K \cap \omega_K$.
Thus, from the viewpoint of $B^3_{N}$ $\partial (\bar{\delta})$ respects
our earlier scheme of when traveling from white to black (resp. black to white) regions having the over-strand to the left (resp. right).
One additional normal position condition we can require by general position arguments is that any $\partial (\bar{\delta})$ intersects any intersection
arc of $\beta_K \cap \omega_K$ at most once.  Now, observe that if we have a curve $\partial (\bar{\delta})$ occurring on one side
of an over-strand at a crossing then we must see a curve $\partial (\bar{\delta^\prime})$ on the another side of the same over-strand.
Finally, observe that if $\partial (\bar{\delta}) = \partial (\bar{\delta^\prime})$ then $\Sigma$ will contain a s.c.c. that is a meridian of $K$.
\begin{figure}[h]
    \includegraphics[height=2.2cm]{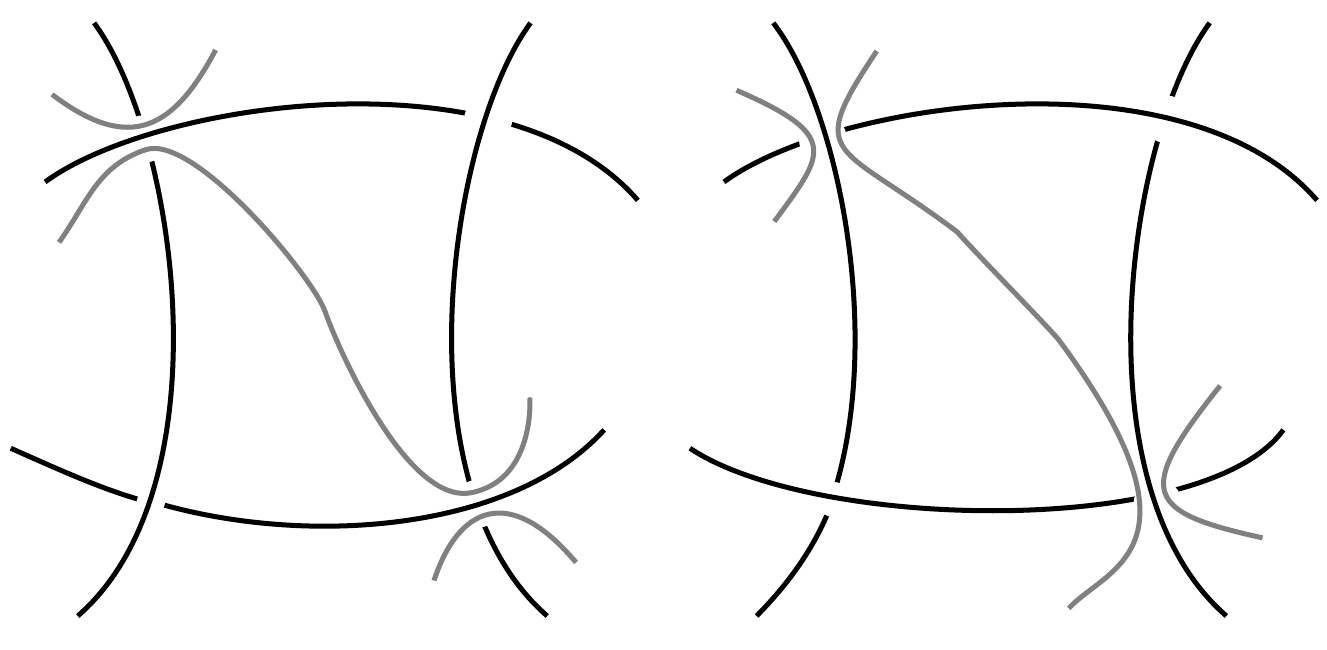}
  \caption{The view of $\partial (\bar{\delta})$ inside $B^3_{N}$.  Traveling along $\partial (\bar{\delta})$, we find a crossing to left/right then right/left.}
  \label{nesting}
\end{figure}
Once placed in the above described normal position, it is readily observed by a nesting argument that $\Sigma$ must contain
a meridian curve  of $K$.  If we walk along the boundary curve of any dome, starting in a white region of $\omega_K$ and passing into
a black region of $\beta_K$, our scheme dictates that the over-strand of a crossing be to our left; then traveling through a black region into a white
region our scheme says the over-strand of a crossing is to our right.  (See figure \ref{nesting}.)  After strolling along all of our boundary curve,
if we have not walked on the two sides of the some crossing then there are at least two
additional curves that are boundaries curves of differing domes, one to the left side of our first crossing and one to the right side of
our second crossing.  Walking along either one of these boundary curves and iterating this procedure we will either find a new boundary of a dome
that we have not strolled along or we will find a boundary that passes on both sides of the over-strand of a crossing, thus yielding a
meridian of $K$.  By the compactness of $\Sigma$ there are only finitely many domes in $B^3_N$ so the latter must occur.

\begin{figure}[h]
    \includegraphics[width=4.8cm]{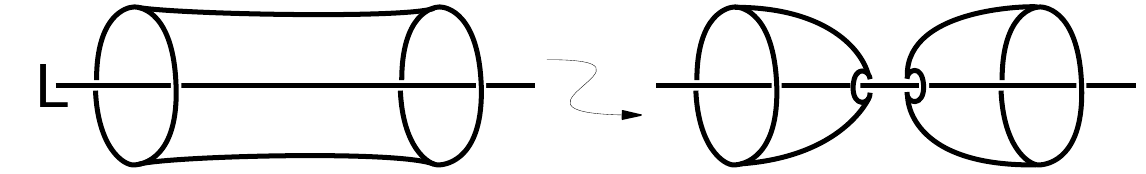}
  \caption{Pairwise compressing $\Sigma$ along a meridian disc.}
  \label{p-compress}
\end{figure}

The implication of an incompressible $\Sigma$ always containing a meridian curve are significant.
Once a meridian curve has been found we can trade $\Sigma$ in for a new surface punctured by $K$, $\Sigma^\prime$,
by compressing $\Sigma$ along a once puncture disc that the meridian curve bounds.  (See figure \ref{p-compress}.)  Repeating this meridian-surgery whenever
a new meridian curve is found, the study of closed incompressible surfaces in the complement of alternating knots becomes the study
of incompressible meridionally or \emph{pairwise incompressible} surfaces---if a curve in $\Sigma$ bounds a disc (resp. once punctured disc) in $S^3 -K$ then
it bounds a disc (resp. once punctured disc) in $\Sigma$.
Finally, the above normal position and nesting argument can be enhanced so as to handle the intersection of an incompressible
pairwise incompressible surface $\Sigma$.  Specifically, if $K$ is an alternating composite knot then there exists a twice punctured $2$-sphere
the intersects $\beta_K$ and $\omega_K$ in single arcs whose union forms a circle that can be thought of as illustrating the composite nature of
$D_K$.

Thus, we have the result that an alternating knot is a composite knot if and only if
its diagram is also composite.  When combined with William P. Thurston's characterization of hyperbolic knots in $S^3$, an immediate consequence is
that all alternating knots coming from diagrams that are neither composite nor $(2, q)$-torus knots are hyperbolic knots.

\section*{Hyperbolic geometry and alternating knots}
\label{HG}
In 1975 Robert Riley produced the first example of a hyperbolic knot, the figure eight, by showing that its fundamental group had a
faithful discrete representation into $PSL( 2 , \mathbb{C})$, the group of orientation preserving isometries of hyperbolic $3$-space.
Inspired by Riley's pioneering work (see page 360 of reference 2), Thurston proved that all knots that are not torus knots and not satellite knots are hyperbolic---their complement has a complete hyperbolic structure.  Recalling that a satellite knot is one where there is an incompressible torus that is not a peripheral torus,
we now apply the fact that for alternating knots such a torus has a meridian curve.  Then pairwise compressing such a torus will produce an annulus
illustrating that the knot is a composite knot.  As mentined previously, we then conclude that all prime alternating knots that are knot $(2,q)$-torus knots are hyperbolic.

Similar to Riley, Thurston was able to structure the hyperbolic structures for knot complements.
However, he used a geometric approach to see their hyperbolic structures.
Specifically, the complement of every hyperbolic knot can be decomposed in a canonical way into ideal tetrahedra, that is,
convex tetrahedra in hyperbolic $3$-space such that all the vertices are lying on the sphere at infinity.

\begin{figure}[h]
    \includegraphics[height=3cm]{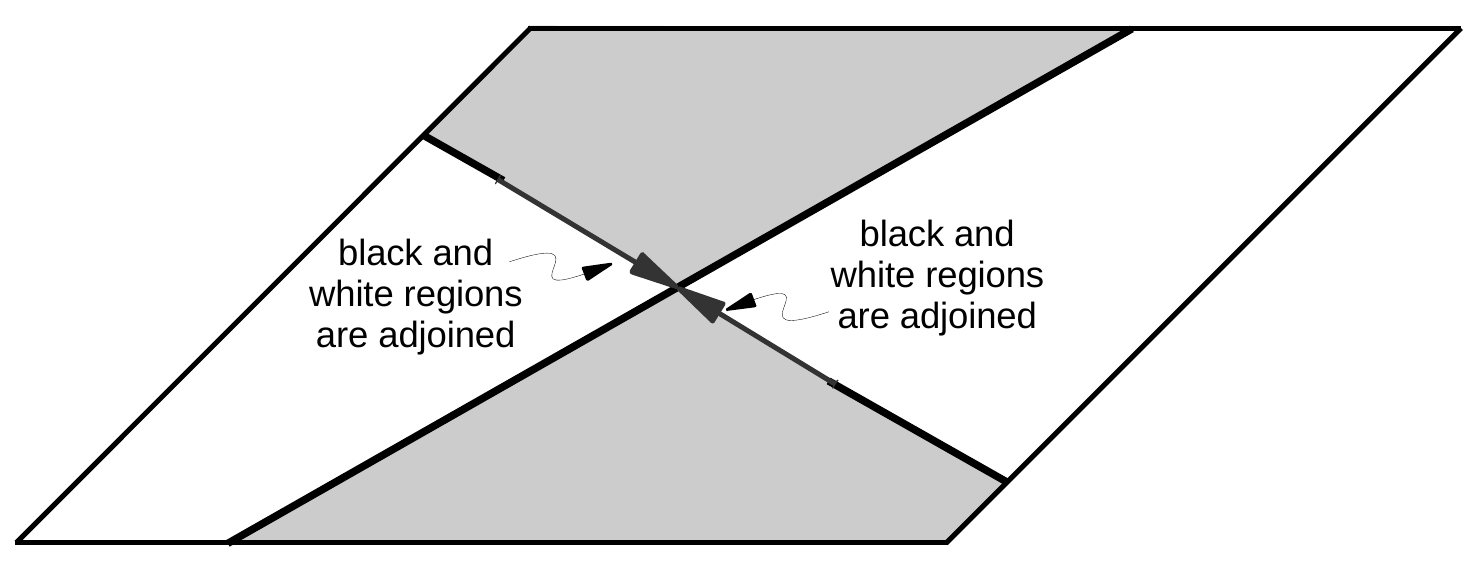}
  \caption{From the viewpoint of $B^3_N$, black and white regions are adjoined along an (oriented) edge).}
  \label{Polyhedron}
\end{figure}

The most famous example of such a tetrahedron decomposition is Thurston's decomposition of the figure eight complement into two regular ideal tetrahedra.
of our previous two open $3$-balls, $B^3_N$ and $B^3_S$.  
Referring to figure \ref{Polyhedron}, if we are viewing the diagram from inside $B^3_N$ and we try to look into $B^3_N$, our view is blocked
by the non-opaque black and white disc-regions of $\beta_{K}$ and $\omega_{K}$.
Moreover, when we look at an arc of
$\beta_{K} \cap \omega_{K}$, from the perspective of $B^3_N$ it does not appear that there are two black (white) regions adjoined be a half-twisted
band.  Instead it appears that the black regions are adjoined to the white regions in a manner consistent with our original scheme
of going from a white region and black region and having the over-strand to our left.  This adjoining scheme for assembling the boundary
of $B^3_N$, and similarly $B^3_S$ from the disc-regions of our colored surfaces
makes sense once one realizes that one's vision in $B^3_N$ gives us information for only the over-strand
at each crossing since we see it in its entirety.  To obtain the information for the under-strand of a crossing---to see it in its entirety---we need
to view it in $B^3_S$ where it appears as an over-strand of a crossing.
\begin{figure}[h]
{\includegraphics[width=6.5cm ]{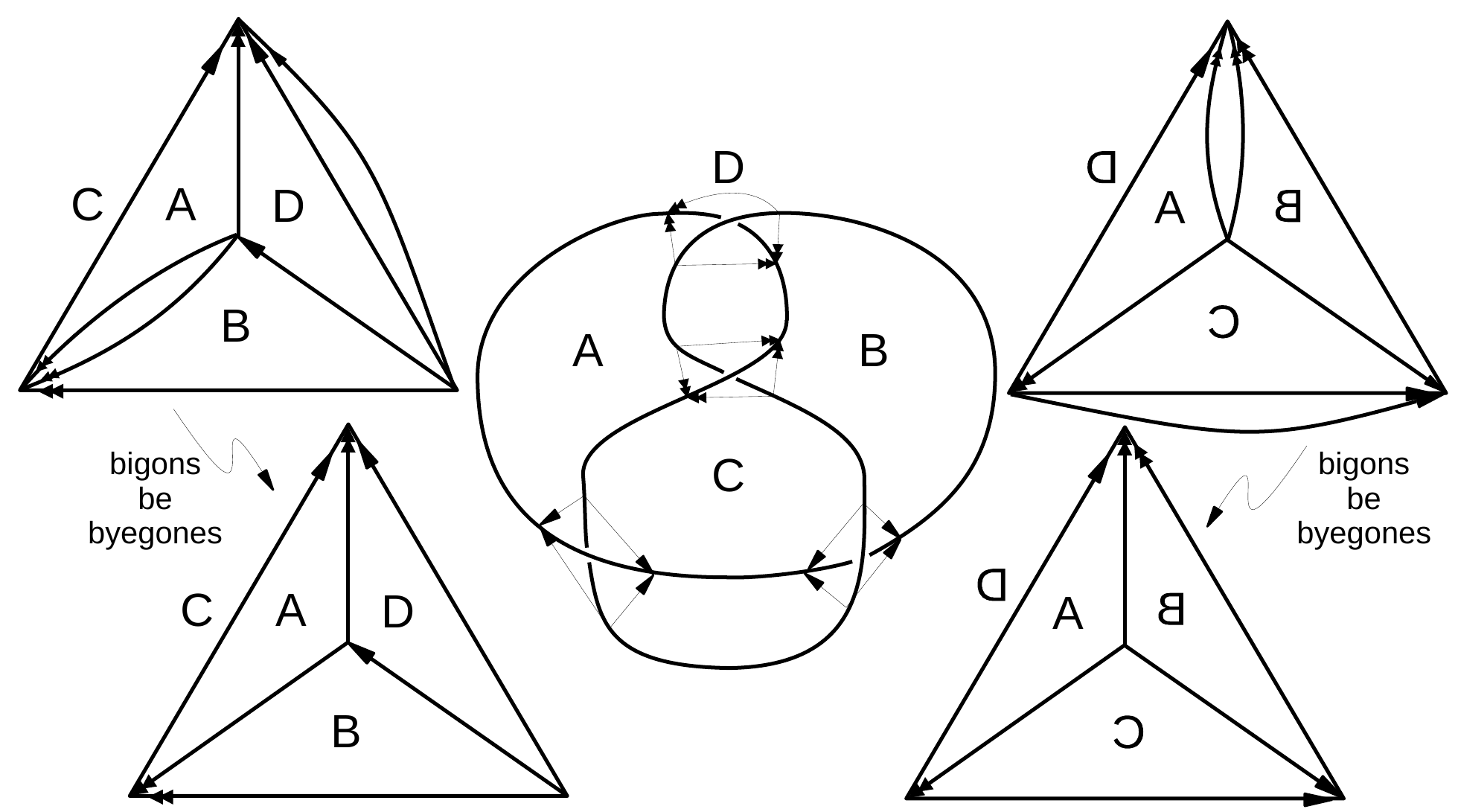}}
\caption{The upper left (resp. right) configurations are the view of $S_N$ (resp. $S_S$) from inside $B^3_N$ (resp. $B^3_S$).}
\label{figure-8}
\end{figure}
In figure \ref{figure-8} we show the inside appearance of the $2$-sphere boundaries, $S_N= \partial B^3_N$ and $S_S = \partial B^3_S$,
coming from this described scheme for assembling of the disc-regions of our colored surface.  Notice that $S_N$ and $S_S$ each
has four vertices---one for each crossing.  When the identification of commonly labeled disc-regions of $S_N$ and $S_S$ is made
these $8$ vertices will become one vertex and the resulting space will correspond to the space obtained by taking the complement
an open tubular neighborhood of the figure eight
and coning its peripheral torus to a point.  Thus, when we delete this single vertex we have the knot complement.
When we place $B^3_N$ and $B^3_S$ in hyperbolic $3$-space, $\mathbb{H}^3$, so that the vertices of $S_N$ and $S_S$ are
in $S^\infty$, the sphere at infinity---they are ideal polyhedra---we have achieved this vertex deletion.  A complete hyperbolic structure for the knot complement will come from a tessellation of $\mathbb{H}^3$ by isometric copies of our two ideal polyhedra.
Now, one feature of an ideal polyhedron is that each edge is the unique geodesic line between two points in $S^\infty$.
Currently since some edges are boundaries of bigon disc-regions in $S_N$ and $S_S$ we need to collapse them to achieve the needed unique edge---we let
bigons be byegones.  For the figure eight, once the four bigon disc-regions are collapsed on $S_N$ and $S_S$, $B^3_N$ and $B^3_S$ become
Thurston's two tetrahedra decomposition.

The above assembly scheme for producing a ideal polyhedron works for any alternating knot that is not a $(2,q)$-torus or composite knot with only one caveat.  If
a disc-region is not a triangle, we need to insert additional edges until it is a union of triangles.

\section*{Establishing conjecture T1}
\label{T1}
In 1985 Vaughan Jones announce the discovery of his polynomial, $V_K(t)$, a Laurent polynomial that is an invariant of
oriented knots (and links).  His discovery emerged from a study of finite dimensional von Neumann algebras.
A number of topologists subsequently reframed the Jones polynomial into more knot theory user friendly settings.  Then in 1987 using
these reframings of the Jones polynomial, Louis Kauffman,
Kunio Murasugi, and Morwen B. Thistlethwaite each independently published proofs of the conjecture T1, \emph{a reduced alternating diagram is a minimal
crossing diagram}.  Additionally, both Mursugi and Thiethlewaite showed that the writhe was an invariant of a
reduced alternating diagram, establishing conjecture T2.

\begin{figure}[h]
    \includegraphics[height=1.6cm]{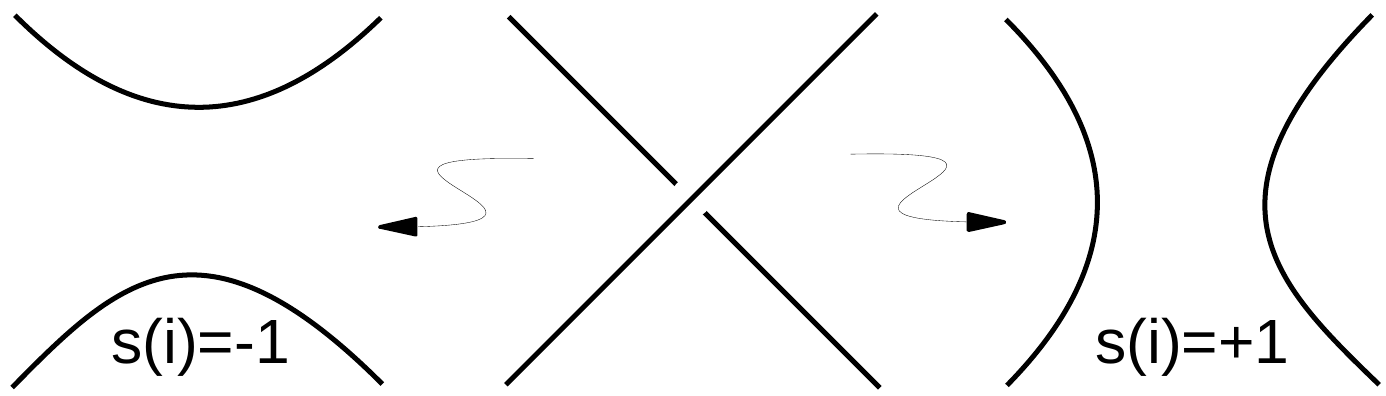}
  \caption{Kauffman's states of $D$.}
  \label{States}
\end{figure}

Kauffman's approach using his bracket invariant and state model equation is particularly accessible.  Moreover, his approach illustrates the importance
of our two colored surfaces.  We consider an $n$-crossing knot diagram, $D$, where we have given each crossing a label $\{1, \cdots , n \}$.
A \emph{state} of $D$ is a choice of smoothly resolving every crossing.  The $i^{\rm th}$ crossing can be resolve either positively or negatively
as shown in figure \ref{States}.

Thus, a state of $D$ also corresponds to a function $s : \{ 1 , \cdots , n \} \rightarrow \{ -1 , +1\}$ where $s(i) = +1 \ ({\rm resp.} \ -1)$
if we smooth the $i^{\rm th}$ crossing in a positive (resp. negative) manner.  Applying a state $s$ to $D$ we obtain a new diagram, $sD$, that
is just a collection of disjoint simple closed curves---there are no crossings.  We will use $|sD|$ to mean the number of curves in the collection.

Of immediate interest is when $D$ is a reduced alternating diagram and the state is constant, $s_+ \equiv +1$.  A constant $+1$ state can be
thought of as cutting every half-twisted band of $\omega_K$.  Notice then $|s_+ D|$ is just equal to the number of disc-regions of $\omega_K$.
Similarly, for $|s_- D|$ is a count of the number of disc-regions of $\beta_K$.  Moreover, using the classical Euler characteristic relation for the sphere,
$ ({\rm vertices} - {\rm edges} + {\rm faces}) = 2$, we can derive the equality $|s_+ D| + |s_- D| = n +2$.

We now come to the \emph{Kauffman bracket} of $D$ which is a Laurent polynomial, $\langle D \rangle $, with an indeterminate variable $A$.
Initially Kauffman axiomatically define his bracket in terms of skein relations.  He then showed it was an invariant of \emph{regular isotopies},
isotopies that corresponded to a sequence of only Reidemeister type II and III moves.  Later, Kauffman derived the following
state sum model for his bracket
\begin{equation}
\label{KB}
\langle D \rangle = \sum_s ( A^{\sum_{i=1}^n s(i)} ( -A^{-2} - A^2 )^{|sD| - 1 } ).
\end{equation}
Although the Kauffman bracket is only a regular isotopy invariant, when we insert a term that takes the writhe, $w(D)$, into account
we do obtain an oriented knot invariant, the \emph{Kauffman polynomial}:
\begin{equation}
\label{KP}
F_K(A) = (-A)^{-3w(D_K)} \langle D_K \rangle.
\end{equation}
As an aside we mention that the Jones polynomial can be obtained from the Kauffman polynomial by a variable substantiation, $F_K( t^{\frac{1}{4}}) = V_K(t)$.

Now since $F_K(A)$ is an oriented knot invariant, the difference between the highest exponent power, $M(F_K)$, and lowest exponent power, $m(F_K)$, of $A$ in $F_K(A)$
is also an invariant.  But, since the $(-A)^{-3w(D_K)}$ factor in equation (\ref{KP}) is common to both $M(F_K)$ and $m(F_K)$ we see that
$M(F_K) - m(F_K) = M(\langle D_K \rangle) - m(\langle D_K \rangle)$, the difference of the highest and lowest power of $A$ in $\langle D_K \rangle$.
Now referring back to equation (\ref{KB}) with a little arguing one can obtain that $M(\langle D_K \rangle) = n + 2|s_+ D| -2$ and
$m(\langle D_K \rangle) = -n -2|s_- D| +2$.  Then:
$$ M(\langle D_K \rangle) - m(\langle D_K \rangle) 
= 2n + 2(|s_+ D| + |s_- D|) -2.$$
Finally, recalling $|s_+ D| + |s_- D| = n +2$ we obtain $M(\langle D_K \rangle) - m(\langle D_K \rangle) = 4n$.

A variation of the above line of reasoning can establish that for an $n$-crossing non-alternating prime diagram, $D$,
we have $ M(\langle D \rangle) - m(\langle D \rangle) < 4n $.  But, since $ M(\langle D \rangle) - m(\langle D \rangle)$ is an invariant
we conclude that a reduced alternating diagram has minimal number of crossings.

\section*{Establishing conjectures T2 \& T3}
\label{T3}
Following the successful verification of conjecture T1 in 1987, Thistlethwaite and William W Menasco announced in 1991 a proof verifying conjecture T3,
the Tait flyping conjecture.  As mentioned before T1 and T3 together imply T2.  Their proof of T3 had two components.  First, it utilized
the new understanding of the interplay
between the algebraic invariants coming out of the Jones polynomial revolution and the crossing number of knot diagrams.
Second, it utilized enhancements of the elementary geometric technology (previously described)  in order to analyze
the intersection pattern of a closed incompressible surfaces with the colored surfaces, $\beta \cup \omega$, of an alternating diagram.

Two key lemmas of this enhanced geometric technology are worth mentioning.  First, for a reduced alternating diagram the associated colored
surfaces, $\beta$ and $\omega$, are themselves incompressible.  Second, for a reduced alternating diagram neither $\beta$ nor $\omega$ contain
a meridional simple closed curve.

The overall strategy of the proof is fairly straight forward.  Given two reduced alternating diagrams, $D_1$ and $D_2$, representing the same
knot type, take the associated two sets of incompressible colored surfaces, $\beta_1 \cup \omega_1$ and $\beta_2 \cup \omega_2$, and
consider how they will intersect each other.  Using nesting arguments that are similar in spirit to the one for closed surfaces, one can establish
that intersections are arcs occurring away from $\beta_1 \cap \omega_1$ and $\beta_2 \cap \omega_2$---the portion of our colored
surfaces that come from the half-twisted bands.  The argument proceeds by trying to ``line up'' crossings of $D_1$ with crossings of $D_2$.
An extensive analysis of these arc intersections then leads to the conclusion
that if $D_1$ and $D_2$ are not isotopic diagrams, there exists a flype.  After performing a sufficient number of the flypes
we have more crossings lining up.  When all the crossings of $D_1$ line up with those of $D_2$ we have isotopic diagrams.

With the establishment of T3 we have the structure for a complete classification of alternating knots.

\section*{What is an alternating knot?}
\label{answer}

In 2017 independently Joshua E. Greene and Joshua Howie published differing answers to Fox's fundamental question.
Howie's answer starts in the setting of knots in $S^3$ and supposes that such a knot, $K$, admits two spanning surfaces, $\Sigma$ and $\Sigma^\prime$.
Neither surface need be orientable.  Let $i(\partial \Sigma , \partial \Sigma^\prime )$ be the intersection number of the the two curves
obtained by intersecting $\Sigma \cup \Sigma^\prime$ with a peripheral torus.  Then Howie proves that
$K$ is an $n$-crossing alternating knot in $S^3$ if and only if the following Euler characteristic equation is satisfied:
$$\chi ( \Sigma ) + \chi ( \Sigma^\prime ) + \frac{1}{2} i(\partial \Sigma , \partial \Sigma^\prime ) = 2.$$
The ``only if'' direction of this statement makes sense once one realizes that $ \frac{1}{2} i(\partial \beta , \partial \omega ) = n$, since each arc of
$\beta \cap \omega$ accounts for two intersections on the peripheral torus, and $\chi ( \Sigma ) + \chi ( \Sigma^\prime ) = (|s_+ D| - n) + (|s_- D| - n)$.

Greene's answer starts in a more general setting, knots in a $\mathbb{Z}_2$-homology $3$-spheres, $M^3$.  He also supposes that one
has two spanning surfaces, $\Sigma$ and $\Sigma^\prime$, of a knot, $K \subset M^3$.  One then considers the Gordon-Litherland pairing form
of the $1^{\rm st}$-homology for each surface:
$$ \mathcal{F}_\Sigma : H_1 (\Sigma) \times H_1 (\Sigma) \rightarrow \mathbb{Z} \ {\rm and} \ 
 \mathcal{F}_{\Sigma^\prime}: H_1 (\Sigma^\prime) \times H_1 (\Sigma^\prime) \rightarrow \mathbb{Z}.$$
\bigskip
If $\mathcal{F}_\Sigma$ is a positive definite form and $\mathcal{F}_{\Sigma^\prime}$ is a negative definite form then Greene's result
states that $M^3 = S^3$ and $K$ is an alternating knot.
Moreover, $\Sigma$ and $\Sigma^\prime$ are our colored surfaces---our two heroes---for an alternating diagram, $D_K$.

\newpage
\begin{center}
{\bf References and further readings}
\end{center}

{\bf Introduction.}
\begin{itemize}
\item[1.] P. G. Tait, On knots I, II, III, Scientific Papers, Cambridge University Press, 1898-1900. Including: Trans. R. Soc. Edin., 28, 1877, 35-79. Reprinted by Amphion Press, Washington D.C., 1993.
\item[2.] T. P Kirkman, The enumeration, description and construction of knots of fewer than ten crossings, Trans. Roy. Soc. Edinburgh 32 (1885), 281-309.
\item[3.] T. P Kirkman, The 364 unifilar knots of ten crossings enumerated and defined, Trans. Roy. Soc. Edinburgh 32 (1885), 483-506.
\item[4.] C. N. Little, On knots, with a census of order ten, Trans. Connecticut Acad. Sci. 18 (1885), 374-378.
\item[5.] C. N. Little, Non alternate ± knots of orders eight and nine, Trans. Roy. Soc. Edinburgh 35 (1889), 663-664.
\item[6.] Jim Hoste, Morwen B. Thistlethwaite and Jeffrey Weeks, The first 1,701,936 knots, Math. Intelligencer 20 (4) (1998), 33-48
\end{itemize}

{\bf The Tait conjectures \& Surfaces and alternating knots.}
\begin{itemize}
\item[1.] William W. Menasco, Closed incompressible surfaces in alternating knot and link complements, Topology 23 (1) (1984), 37-44.
\item[2.] William W. Menasco, Determining incompressibility of surfaces in alternating knot and link complements. Pacific J. Math. 117 (1985), no. 2, 353-370.
\item[3.]  William W. Menasco and M. B Thistlethwaite, Surfaces with boundary in alternating knot exteriors. J. Reine Angew. Math. 426 (1992), 47-65.
\end{itemize}

{\bf Hyperbolic geometry and alternating knots.}
\begin{itemize}
\item[1.] Robert Riley, Discrete parabolic representations of link groups, Mathematika, 22 (2) (1975), pp. 141-150.
\item[2.] William P. Thurston, Three dimensional manifolds, Kleinian groups and hyperbolic geometry, Bull. Amer. Math. Soc. (NS), 6 (1982), pp. 357-381
\item[3.] William W. Menasco, Polyhedra representation of link complements. Low-dimensional topology (San Francisco, Calif., 1981), 305?325, Contemp. Math., 20, Amer. Math. Soc., Providence, RI, 1983.
\item[4.]  Marc Lackenby, The volume of hyperbolic alternating link complements. With an appendix by Ian Agol and Dylan Thurston. Proc. London Math. Soc. (3) 88 (2004), no. 1, 204?224.
\end{itemize}

{\bf Establishing conjecture T1.}
\begin{itemize}
\item[1.] Louis H. Kauffman, State models and the Jones polynomial, Topology 26 (1987), no. 3, 395-407.
\item[2.] Kunio Murasugi, Jones polynomials and classical conjectures in knot theory, Topology 26 (1987), no. 2, 187-194.
\item[3.] Morwen B. Thistlethwaite, A spanning tree expansion of the Jones polynomial, Topology 26 (1987), no. 3, 297-309.
\item[4.] Vaughan Jones, A Polynomial Invariant for Knots via von Neumann Algebras. Bull. Am. Math. Soc. 12, 103-111, 1985.
\item[5.] Vaughan Jones, "Hecke Algebra Representations of Braid Groups and Link Polynomials." Ann. Math. 126, 335-388, 1987.
\item[5.] W.B.R. Lickorish, An introduction to knot theory, Graduate Texts in Mathematics Series, Springer Verlag, 1997.
\end{itemize}

{\bf Establishing conjectures T2 \& T3.}
\begin{itemize}
\item[1.] William W. Menasco and M. B Thistlethwaite, The classification of alternating links, Ann. Math. 138 (1993), 113-171.
\item[2.]  William W. Menasco and M. B Thistlethwaite, The Tait flyping conjecture. Bull. Amer. Math. Soc. (N.S.) 25 (1991), no. 2, 403-412.
\end{itemize}

{\bf What is an alternating knot?}
\begin{itemize}
\item[1.] Joshua A. Howie, A characterisation of alternating knot exteriors. Geom. Topol. 21 (2017), no. 4, 2353-2371.
\item[2.] Joshua E. Greene, Alternating links and definite surfaces. With an appendix by András Juhász and Marc Lackenby. Duke Math. J. 166 (2017), no. 11, 2133-2151.
\end{itemize}

\end{document}